\documentclass[reqno,preprint,11pt]{elsarticle}
\usepackage[active]{srcltx}
\usepackage{bbm}
\usepackage{tipa}
\usepackage{amsthm}
\usepackage{amssymb}
\usepackage{stackrel}
\usepackage{mathrsfs}
\usepackage{amsmath,amscd}
\usepackage{lscape}
\usepackage{longtable}
\usepackage{xltabular} 
\usepackage{ragged2e}
\usepackage{booktabs}
\usepackage{rotating}
\usepackage[active]{srcltx}
\usepackage{amsfonts}
\usepackage[centerlast]{caption}
\usepackage{xcolor}
\definecolor{myred}{RGB}{200, 0, 32}
\definecolor{backgrey}{HTML}{E6E6E6} 
\usepackage{latexsym}
\usepackage{amsmath,amsthm}
\usepackage{amssymb,amsmath}
\usepackage[all]{xy}
\usepackage{graphicx}
\usepackage[dvips]{geometry}

\newtheorem{dl}{Theorem}[section]
\newtheorem{tl}[dl]{Corollary}
\newtheorem{yl}[dl]{Lemma}

\newtheorem{xz}[dl]{Proposition}

\numberwithin{equation}{section}
\newproof{pot1}{\bf Proof of Theorem \ref{GenAI}}
\newproof{pot2}{\bf Proof of Theorem \ref{GenAII}}
\newproof{pot3}{\bf Proof of Theorem \ref{andrewsthree}}
\newproof{pot4}{\bf Proof of Theorem \ref{warnaar}}
\newproof{pot5}{\bf Proof of Theorem \ref{war-one}}
\def\qed{\hfill \rule{4pt}{7pt}}
\def\pf{\noindent {\sl Proof.~}}
\newcommand{\poq}[2]{(#1;q)_{#2}}

\newcommand{\poqq}[2]{(#1;q^2)_{#2}}

\usepackage[colorlinks=true,linkcolor=green,filecolor=brown,citecolor=green]{hyperref}
\usepackage{hyperref}
\usepackage{epsf}

\allowdisplaybreaks
\title{Some new results on Andrews' and Warnaar's $q$-identities}
\author{Qi Chen~~\fnref{fn22,fn3}}
\address[fn22]{Department of Mathematics, Soochow University, SuZhou 215006, P.R.China}
\fntext[fn3]{E-mail address: qchen2025@stu.suda.edu.cn}
\begin{document}
\begin{abstract}
In this  paper, by the technique of inverse relations and comparing coefficients, we  establish some generalized forms of Andrews' $q$-series identity and two new Bailey pairs and $q$-identities closely related to Andrews-Warnaar's sum identity for partial theta functions. 
\end{abstract}

\maketitle
\section{Introduction}
In the theory of $q$-series, there is a well-known identity given by G.E. Andrews.
\begin{dl}[Andrews' identity: {\rm \cite[(II. 11)]{10}}]\label{AI}
\begin{align*}
     \sum_{n=0}^\infty\frac{\poq{a,b}{n}}{\poq{q}{n}\poqq{abq}{n}}q^{\frac{n(n+1)}{2}}=\frac{\poq{-q}{\infty}\poqq{aq,bq}{\infty}}{\poqq{abq}{\infty}}.
\end{align*}
\end{dl}
A quick glance on this identity shows that
 the series expansion of the factors 
  \begin{align}
        \frac{\poqq{ax}{\infty}}{\poq{x}{\infty}}\quad\mbox{and}\quad  \frac{\poq{ax}{\infty}}{\poqq{x}{\infty}}\label{q1/q}
          \end{align}
          seem necessary but to have been ignored, merely comparing with the $q$-binomial theorem \cite[(II. 3)]{10}
    \begin{align}
        \frac{\poq{ax}{\infty}}{\poq{x}{\infty}}&=\sum_{n=0}^\infty\frac{\poq{a}{n}}{\poq{q}{n}}x^n.\label{q2/q}
    \end{align}
The latter has been frequently used in study of $q$-series. Therefore, it is natural to consider the former's expansion as power series of $x$. We believe such expansions also play a similar role as the $q$-binomial theorem in the study of $q$-series.  

As an initial positive answer to this question, in this paper we will establish the following $q$-identities,  each of which can be regarded as generalizations of Theorem \ref{AI}.
\begin{dl}[Generalization of Andrews' identity-I] \label{GenAI} For any complex numbers $a, b,$ and $c$, there  holds
\begin{align}
     \sum_{n=0}^\infty\frac{\poq{a,b}{n}}{\poq{q}{n}\poqq{abq}{n}}q^{\frac{n(n+1)}{2}}{c^n}=\frac{\poq{-cq,a,b}{\infty}}{\poqq{abq}{\infty}}\sum_{n=0}^\infty\frac{h_n(a,b|q^2)}{\poq{q,-cq}{n}},\label{GenAID}
\end{align}
where the Rogers-Szeg{\"o} polynomials are defined by
\begin{align}
    h_n(a,b|q):=\sum_{k=0}^n\genfrac{[}{]}{0pt}{}{n}{k}_{q}a^kb^{n-k}.\label{chenqi}
\end{align}
\end{dl}

\begin{dl}[Generalization of Andrews' identity-II] \label{GenAII}Let $h_n(a,b|q)$ be defined by \eqref{chenqi}. We have
\begin{align*}
   &\sum_{n=0}^\infty\frac{\poq{a,b}{n}}{\poq{q}{n}\poqq{ab}{n}}q^{\frac{n(n+1)}{2}}{c^n}=\frac{\poq{a,b,-cq}{\infty}}{\poqq{ab}{\infty}}\sum_{n=0}^\infty\frac{h_n(a,bq|q^2)+h_n(aq,b|q^2)}{\poq{q,-cq}{n}(1+q^n)}.
\end{align*}
\end{dl}
Similar to but slightly different from these two conclusions is the following $q$-identity.
\begin{dl}
\label{andrewsthree}
\begin{align}
 &  \sum_{n=0}^\infty \frac{\poq{a,b}{n}}{\poq{q,(ab)^{1/2}}{n}\poq{-(ab)^{1/2}}{n+1}}q^{\frac{n(n+1)}{2}}c^n\label{andrews-twotwo}\\
 &=\frac{\poq{a,b,-cq}{\infty}}{\poqq{ab}{\infty}}\sum_{n=0}^\infty\bigg(\frac{a^{1/2}}{a^{1/2}+b^{1/2}} \frac{h_n(a,bq|q^2)}{\poq{q,-cq}{n}}+\frac{b^{1/2}}{a^{1/2}+b^{1/2}}\frac{h_n(aq,b|q^2)}{\poq{q,-cq}{n}}\bigg).\nonumber
\end{align}
\end{dl}
It is worth mentioning that the case $c=1$ of Theorem \ref{andrewsthree} is very similar to \cite[Thm. 1]{222}, 
the latter is used to  generating functions for the partitions of asymmetric residue classes. The reader is referred to loc. cit for more details. 

Here and in what follows, we will  adopt the standard notation and terminology for basic hypergeometric series in the book \cite{10}. Let $|q|<1$ and the $q$-shifted
factorials be given by
\begin{align*}
(x;q)_{\infty}=\prod_{n=0}^{\infty}(1-xq^n),\quad (x;q)_n=\frac{(x;q)_{\infty}}{(xq^n;q)_{\infty}}.
\end{align*}
For convenience, we will adopt the following notation for multiple $q$-shifted factorial
\begin{align*}
(x_1,x_2,\cdots,x_m;q)_{n}=(x_1;q)_{n}(x_2;q)_{n}\cdots(x_m;q)_{n}.
\end{align*}
The basic hypergeometric series with the base $q$ and the argument $x$ is defined by
\begin{align}
_{r}\phi _{s}\left(\begin{matrix}a_1,&a_2, \cdots, &a_r \\ b_1,&b_2, \cdots, &b_s \end{matrix}; q, x \right)
=\sum_{n=0}^{\infty}\frac{(a_1,a_2,\cdots,a_r;q)_{n}}{(q,b_1,b_2,\cdots,b_s;q)_{n}}
(\tau_1(n))^{s+1-r}x^n, \label{(1.8)}
\end{align}
where  $\tau_r(n)=(-1)^nq^{rn(n-1)/2}$, $\{b_i\}$ are complex numbers such that none of the denominators in \eqref{(1.8)} vanishes, which we never mention as conditions in what follows unless otherwise stated.  As a custom,
the $q$-binomial coefficient $\genfrac{[}{]}{0pt}{}{n}{k}_{q}$ is defined by
 \begin{align*}
\genfrac{[}{]}{0pt}{}{n}{k}_{q}:=
 \displaystyle  \frac{\poq{q}{n}}{\poq{q}{k}\poq{q}{n-k}}.
\end{align*}

Our paper is organized as follows. In the next section we will show some preliminaries and then give the full proofs for Theorems \ref{GenAI}, \ref{GenAII}, and \ref{andrewsthree}. After that, some concrete $q$-identities will be discussed.  In Section \ref{sec3}, we investigate Warnaar's partial theta identities in the same vein, two  new Bailey pairs and some $q$-identities are presented.
 \section{Proofs of lemmas and theorems}\label{sec2}
  \subsection{Preliminary lemmas}  
    As a positive answer to the expansion question of \eqref{q1/q} and  one of the necessary preliminary results, we now  show
\begin{yl}\label{main}Suppose that
     \begin{align}
        \frac{\poqq{ax}{\infty}}{\poq{x}{\infty}}&=\sum_{n=0}^\infty\lambda_n(a)x^n.\label{expanques}
    \end{align}
    Then the coefficients
    \begin{align}
        \lambda_n(a)&={}_2\phi_1\left(\begin{array}{cc}
        q^{-n}, & -q^{-n} \\
         & 0
    \end{array};q,q^2/a\right)\frac{(-a)^nq^{n(n-1)}}{\poqq{q^2}{n}}.\label{coeff}
    \end{align}
\end{yl}
\pf  It is obvious that \eqref{expanques} is equivalent to
\begin{align*}
       \frac{1}{\poq{x}{\infty}}=\frac{1}{\poqq{ax}{\infty}}\sum_{n=0}^\infty\lambda_n(a) x^n.
     \end{align*}
    Namely
\begin{align*}
 \sum_{n=0}^{\infty}\frac{x^n}{\poq{q}{n}}=\sum_{k=0}^\infty \frac{a^kx^k}{\poqq{q^2}{k}}\sum_{n=0}^\infty\lambda_n(a) x^n.
 \end{align*}
 By equating the coefficients of $x^n$, we get
 \begin{align*} &\frac{1}{\poq{q}{n}}=\sum_{k=0}^n\frac{a^{n-k}}{\poqq{q^2}{n-k}}\lambda_k(a),
\end{align*}
which is
 \begin{align*} \frac{\poqq{q^2}{n}}{\poq{q}{n}}a^{-n}=\sum_{k=0}^n\genfrac{[}{]}{0pt}{}{n}{k}_{q^2}\poqq{q^2}{k}a^{-k}\lambda_k(a).
\end{align*}
After simplification, it takes the form
\begin{align}
    \alpha_n=\sum_{k=0}^n \genfrac{[}{]}{0pt}{}{n}{k}_{q^2}\beta_k,\label{equation}
\end{align}
where 
$$
    \alpha_n=\poq{-q}{n}/a^{n}, \quad \beta_n=\poqq{q^2}{n}\lambda_n(a)/a^{n}.
$$
Recall that there is a matrix inversion
\begin{align*}
    \left(\genfrac{[}{]}{0pt}{}{n}{k}_{q^2}\right)_{n\geq k\geq 0}^{-1}=\left(\genfrac{[}{]}{0pt}{}{n}{k}_{q^2}\tau_{2}(n-k)\right)_{n\geq k\geq 0}.
\end{align*} By applying this inverse relation, we solve $\lambda_n(a)$ from \eqref{equation} as below
\begin{align*}
    \poqq{q^2}{n}\lambda_n(a) a^{-n}&=\sum_{k=0}^n\genfrac{[}{]}{0pt}{}{n}{k}_{q^2}\tau_{2}(n-k)\poq{-q}{k}a^{-k}\\
    &=\tau_{2}(n)\sum_{k=0}^n\genfrac{[}{]}{0pt}{}{n}{k}_{q^2}\tau_{2}(k)q^{2(k-nk)}\poq{-q}{k}a^{-k}.
\end{align*}
A bit simplification gives
\begin{align*}
    \lambda_n(a)&=\frac{a^n\tau_2(n)}{\poqq{q^2}{n}}\sum_{k=0}^n\frac{\poqq{q^{-2n}}{k}}{\poqq{q^2}{k}}\poq{-q}{k}(q^{2}/a)^{k}\\
    &=\frac{a^n\tau_2(n)}{\poqq{q^2}{n}}{}_2\phi_1\left(\begin{array}{cc}
        q^{-n}, & -q^{-n} \\
         & 0
    \end{array};q,q^2/a\right).
\end{align*}
That is we wanted. \qed

Apart from  Lemma \ref{main}, we also require the following recurrence relation related to $\lambda_n(a)$.
\begin{yl} Define, for integers $r, s$ and $n$, the finite sum
  \begin{align}
   T_{r,n}(s):= \sum_{k=0}^n\frac{\poqq{q^{-2n}}{k}}{\poq{q,q^{1+r-n}}{k}}q^{(2-s)k}.\label{tms}
\end{align}
   Then there is a  recurrence relation of degree two for $\{T_{r,n}(s)\}_{n\geq 0}$
    \begin{align}
      q^{-s}(1-q^{2n+2}) (q^{-s}+q^{r+n})&T_{r,n}(s)\nonumber\\
      +q^{n}(q^n-q^r)(-q^{2n+1}+q^{r+n}+q^{-s}+q^{-s-1})&T_{r,n+1}(s)\nonumber\\+q^{2n} (q^{n}-q^{r}) (q^{n+1}-q^{r})&T_{r,n+2}(s)=0.  \label{tms-rec}
    \end{align}
\end{yl}
\pf  
For convenience, we write $y$ for $q^s$.  Then the recurrence relation \eqref{tms-rec} becomes
\begin{align}
&y^{-1}\left(1-q^{2 n+2}\right)\left(y^{-1}+q^{r+n}\right)  \sum_{k=0}^n \frac{\poqq{q^{-2n}}{k}}{\poq{q,q^{1+r-n}}{k}} q^{2k}y^{-k} \nonumber\\
&+q^n\left(q^n-q^r\right)\left(-q^{1+2 n}+q^{r+n}+(1+q)  q^{-1}   y^{-1}\right)  \sum_{k=0}^{n+1} \frac{\poqq{q^{-2 n-2}}{k}  }{\poq{q,q^{r-n}}{k}} q^{2 k} y^{-k}\nonumber\\
& =- q^{2 n}\left(q^n-q^r\right)\left(q^{n+1}-q^r\right)\sum_{k=0}^{n+2} \frac{\poqq{q^{-2n-4}}{k}}{\poq{q,q^{r-n-1}}{k}} q^{2k} y^{-k}.\label{equivancy}
\end{align}
As such, in order to show \eqref{equivancy} we only need  to compare the coefficients of $y^{-k}$ on both sides. Thus, we have 
\begin{align*}
\boldsymbol\lbrack y^{-k}\boldsymbol\rbrack\text{LHS of \eqref{equivancy}}=&  \left(1-q^{2 n+2}\right) \frac{\left(q^{-2 n} ; q^2\right)_{k-2}}{\left(q, q^{1+r-n} ; q\right)_{k-2}} q^{2(k-2)}(=S_1)\\
&+\left(1-q^{2 n+2}\right) q^{r+n}  \frac{\left(q^{-2 n} ; q^2\right)_{k-1}}{\left(q, q^{1+r-n} ; q^2\right)_{k-1}} q^{2(k-1)}(=S_2) \\
& + q^n\left(q^n-q^r\right)  \left(-q^{1+2 n}+q^{r+n}\right)  \frac{\left(q^{-2 n-2} ; q^2\right)_{k }}{\left(q, q^{r-n }; q)_k\right.} q^{2 k}(=S_3)\\
&+q^n\left(q^n-q^r\right)(1+q)   q^{-1}   \frac{\left(q^{-2n-2} ; q^2\right)_{k-1}}{\left(q,q^{r-n} ; q\right)_{k-1}} q^{2(k-1)}(=S_4).
\end{align*}
Here and in the sequel, the notation $\boldsymbol\lbrack t^n\boldsymbol\rbrack f(t)$ denotes the coefficient of $t^n$ in $f(t)$.
Consequently, it is easy to check 
\begin{align*}
 S_1+S_4&= \frac{\poqq{q^{-2n-4}}{k}}{\poq{q,q^{r-n-1}}{k}}q^{2n+2k-3}\frac{(1-q^{r-n-1})(1-q^{r-n})(1-q^{2k})}{1-q^{-2n-4}},\\
  S_2+S_3&=-\frac{\left(q^{-2 n-2} ; q^2\right)_{k }}{\left(q, q^{r-n }; q)_k\right.} q^{2 k+4n+1}(1-q^{r-n})(1-q^{r-n+k-1}).
\end{align*}
Based on these two results, it follows
\begin{align*}
   \boldsymbol\lbrack y^{-k}\boldsymbol\rbrack\text{LHS of \eqref{equivancy}}&=\frac{\poqq{q^{-2n-4}}{k}}{\poq{q,q^{r-n-1}}{k}} q^{2k}\bigg\{q^{2n-3}\frac{(1-q^{r-n-1})(1-q^{r-n})(1-q^{2k})}{1-q^{-2n-4}}\\
   &\qquad\qquad\qquad -q^{4n+1}\frac{(1-q^{r-n-1})(1-q^{-2n-4+2k})}{1-q^{-2n-4}}(1-q^{r-n})\bigg\}\\
 & =\frac{\poqq{q^{-2n-4}}{k}}{\poq{q,q^{r-n-1}}{k}} q^{2k}(-q^{4n+1})(1-q^{-2n-4})\frac{(1-q^{r-n-1})(1-q^{r-n})}{1-q^{-2n-4}}\\
  &=\boldsymbol\lbrack y^{-k}\boldsymbol\rbrack\mbox{RHS of \eqref{equivancy}}.
\end{align*}
As desired.
\qed

Very interesting is that \eqref{tms-rec} includes  the $q$-Chu-Vandermonde formula \cite[(II. 6)]{10} as a special case.
\begin{yl}\label{lemmaadded}
Let $T_{r,n}(s)$ be given by \eqref{tms}.
 Then the following identities hold.
    \begin{align}
   T_{r,n}(1)&=\frac{\poq{-q^{1+r}}{n}}{\poq{q^{1+r-n}}{n}}(-1)^nq^{-n^2},\label{s1}\\
       T_{r,n}(0)&=\frac{q^n+q^r}{1+q^{r+n}}\frac{\poq{-q^{1+r}}{n}}{\poq{q^{1+r-n}}{n}}(-1)^nq^{-n^2+n}.\label{s2}
           \end{align}
   \end{yl}
  \subsection{The full proofs of the main theorems} 
As is expected, we turn  to present  full proofs for Theorems \ref{GenAI}-\ref{andrewsthree} mainly based on  Lemmas \ref{main}-\ref{lemmaadded}. As the first step, we need to show a master theorem.
\begin{dl}\label{andrews-s} Let $T_{r,n}(s)$ be given by \eqref{tms}. Then
      \begin{align}
     \sum_{n=0}^\infty&\frac{\poq{a,b}{n}q^{\frac{n(n+1)}{2}}{c^n}}{\poq{q}{n}\poqq{abq^s}{n}}=\frac{\poq{a,b,-cq}{\infty}}{\poqq{abq^s}{\infty}}
      \sum_{n\geq k\geq 0}\frac{\tau_2(k)(bq^s)^k}{\poqq{q^2}{k}}\frac{a^{n-k}}{\poq{q}{n-2k}}\frac{T_{n-k,k}(s)}{\poq{-cq}{n}}.\label{andrews-sdl}
 \end{align}
\end{dl}
\pf
 It is clear that from Lemma \ref{main}, it follows
\begin{align*}
   \frac{\poqq{abq^{s+2n}}{\infty}}{\poq{bq^n}{\infty}}=\sum_{M=0}^\infty \frac{(aq^{s+n})^M\tau_2(M)}{\poqq{q^2}{M}}{}_2\phi_1\left(\begin{array}{cc}
         q^{-M}, & -q^{-M} \\
          & 0
     \end{array};q,q^{2-s-n}/a\right)(bq^n)^M.
\end{align*}
Thus
\begin{align*}
    &\sum_{n=0}^\infty\frac{\poq{a,b}{n}q^{\frac{n(n+1)}{2}}{c^n}}{\poq{q}{n}\poqq{abq^s}{n}}
     =\frac{\poq{b}{\infty}}{\poqq{abq^s}{\infty}}\sum_{n=0}^\infty \frac{\poq{a}{n}}{\poq{q}{n}}q^{n(n+1)/2}c^n\\
     &\qquad \times\sum_{M=0}^\infty \frac{(aq^{s+n})^M\tau_2(M)}{\poqq{q^2}{M}}{}_2\phi_1\left(\begin{array}{cc}
         q^{-M}, & -q^{-M} \\
          & 0
     \end{array};q,q^{2-s-n}/a\right)(bq^n)^M.
\end{align*}
Written out in explicit terms, it is
\begin{align*}
     \sum_{n=0}^\infty\frac{\poq{a,b}{n}q^{\frac{n(n+1)}{2}}{c^n}}{\poq{q}{n}\poqq{abq^s}{n}}
         &=\frac{\poq{a,b}{\infty}}{\poqq{abq^s}{\infty}}\sum_{n=0}^\infty \frac{q^{n(n+1)/2}c^n}{\poq{q}{n}}\sum_{M=0}^\infty\frac{q^{(s+n)M}\tau_2(M)}{\poqq{q^2}{M}}\nonumber\\
    &\times\sum_{k=0}^M\frac{\poqq{q^{-2M}}{k}}{\poq{q}{k}}q^{k(2-s-n)}(bq^n)^M  \frac{a^{M-k}}{\poq{aq^n}{\infty}}.
\end{align*}
Rearrange the right-hand infinite series as power series in $a$. We have
\begin{align*}
    \sum_{n=0}^\infty\frac{\poq{a,b}{n}q^{\frac{n(n+1)}{2}}{c^n}}{\poq{q}{n}\poqq{abq^s}{n}}
     &=\frac{\poq{a,b}{\infty}}{\poqq{abq^s}{\infty}}\sum_{n=0}^\infty \frac{q^{n(n+1)/2}c^n}{\poq{q}{n}}\sum_{M=0}^\infty\frac{q^{(s+n)M}\tau_2(M)}{\poqq{q^2}{M}}\nonumber\\
    &\times\sum_{k=0}^M\frac{\poqq{q^{-2M}}{k}}{\poq{q}{k}}q^{k(2-s-n)}(bq^n)^M\sum_{r=0}^\infty \frac{q^{n(r-(M-k))}}{\poq{q}{r-(M-k)}}a^r.
\end{align*}
After simplification and rearrangement  of summations over $n$ and $k$, it becomes
\begin{align*}
    \sum_{n=0}^\infty\frac{\poq{a,b}{n}q^{\frac{n(n+1)}{2}}{c^n}}{\poq{q}{n}\poqq{abq^s}{n}}
    &=\frac{\poq{a,b}{\infty}}{\poqq{abq^s}{\infty}}\sum_{M=0}^\infty\frac{(bq^s)^M\tau_2(M)}{\poqq{q^2}{M}}\sum_{r=0}^\infty \frac{a^r}{\poq{q}{r-M}}\\
     &\times\sum_{k=0}^M\frac{\poqq{q^{-2M}}{k}}{\poq{q,q^{1+r-M}}{k}}q^{k(2-s)}\sum_{n=0}^\infty \frac{q^{n(n-1)/2}c^n}{\poq{q}{n}}q^{n(1+r+M)}.
\end{align*}
Observe that the last summation can be evaluated in closed form by the $q$-binomial theorem \eqref{q2/q}. Thus we get
\begin{align*}
    &\sum_{n=0}^\infty\frac{\poq{a,b}{n}q^{\frac{n(n+1)}{2}}{c^n}}{\poq{q}{n}\poqq{abq^s}{n}}\\
    &=\frac{\poq{a,b,-cq}{\infty}}{\poqq{abq^s}{\infty}}\sum_{M=0}^\infty\frac{(bq^s)^M\tau_2(M)}{\poqq{q^2}{M}}\sum_{r=0}^\infty \frac{a^r}{\poq{q}{r-M}\poq{-cq}{r+M}}
    \sum_{k=0}^M\frac{\poqq{q^{-2M}}{k}q^{k(2-s)}}{\poq{q,q^{1+r-M}}{k}}\\
     &\stackrel{N:=r+M}{===}\frac{\poq{a,b,-cq}{\infty}}{\poqq{abq^s}{\infty}}\sum_{N=0}^\infty \frac{1}{\poq{-cq}{N}}\sum_{M=0}^\infty\frac{(bq^s)^M\tau_2(M)a^{N-M}}{\poqq{q^2}{M}\poq{q}{N-2M}}
     \sum_{k=0}^M\frac{\poqq{q^{-2M}}{k}q^{k(2-s)}}{\poq{q,q^{1+N-2M}}{k}}.
\end{align*}
Reformulating this identity in terms of \eqref{tms}, we get  \eqref{andrews-sdl}.
\qed
 
   First, by utilizing Theorem \ref{andrews-s} and the identity \eqref{s1}, we begin to show Theorem \ref{GenAI}.
\begin{pot1}Observe that the special case  $s=1$ of \eqref{andrews-sdl} yields
\begin{align*}
    \sum_{n=0}^\infty\frac{\poq{a,b}{n}q^{\frac{n(n+1)}{2}}{c^n}}{\poq{q}{n}\poqq{abq}{n}}=\frac{\poq{a,b,-cq}{\infty}}{\poqq{abq}{\infty}}\sum_{n\geq k\geq 0} \frac{\tau_2(k)(bq)^k}{\poqq{q^2}{k}}\frac{a^{n-k}}{\poq{q}{n-2k}}\frac{T_{n-k,k}(1)}{\poq{-cq}{n}}.
\end{align*}
Substituting \eqref{s1} into $T_{n-k,k}(1)$, we can derive
\begin{align*}
    \sum_{n=0}^\infty&\frac{\poq{a,b}{n}q^{\frac{n(n+1)}{2}}{c^n}}{\poq{q}{n}\poqq{abq}{n}}=\frac{\poq{a,b,-cq}{\infty}}{\poqq{abq}{\infty}}
   \sum_{n=0}^\infty\frac{1}{\poq{-cq}{n}}\\ &\times\sum_{k=0}^n\frac{\tau_2(k)(bq)^k}{\poqq{q^2}{k}}\frac{a^{n-k}}{\poq{q}{n-2k}}\frac{\poq{-q^{1+n-k}}{k}}{\poq{q^{1+n-2k}}{k}}(-q^{-k})^k\\
    &=\frac{\poq{a,b,-cq}{\infty}}{\poqq{abq}{\infty}}\sum_{n=0}^\infty\frac{\poq{-q}{n}}{\poq{-cq}{n}}\sum_{k=0}^n\frac{b^ka^{n-k}}{\poqq{q^2}{k}\poqq{q^2}{n-k}}\\
    &=\frac{\poq{a,b,-cq}{\infty}}{\poqq{abq}{\infty}}\sum_{n=0}^\infty\frac{1}{\poq{q,-cq}{n}}\sum_{k=0}^n\genfrac{[}{]}{0pt}{}{n}{k}_{q^2}b^ka^{n-k}.
\end{align*}
Note that the inner sum is the Rogers-Szeg{\"o} polynomials \eqref{chenqi} on the base $q^2$.  The proof is finished.\qed
\end{pot1}

From \eqref{GenAID} of Theorem \ref{GenAI}, we can deduce some interesting $q$-identities. 
For instance,  letting $c=0$ in \eqref{GenAID}, we have a new generating function of the Rogers-Szeg\"{o} polynomials.
\begin{tl}\label{444444}Let $h_n(a,b|q)$ be defined by \eqref{chenqi}. Then we have
    \begin{align}
     \sum_{n=0}^\infty\frac{h_n(a,b|q^2)}{\poq{q}{n}}
     =\frac{\poqq{abq}{\infty}}{\poq{a,b}{\infty}}.\label{newGF}
\end{align}
     \end{tl}
 A combination of this result with Lemma \ref{main} gives rise to a finite $q$-series identity
    \begin{align*}
      (-a)^Mq^{M^2}\sum_{k=0}^M\frac{\poqq{q^{-2M}}{k}}{\poq{q}{k}}(q/a)^k
      =\sum_{k=0}^M\genfrac{[}{]}{0pt}{}{M}{k}_{q}\poq{a}{k}q^{k(k+1)/2},
    \end{align*} 
 which is in agreement with the special case $(b,c,z)=(a,0,-q^{M+1})$ of \cite[(III. 6)]{10}.
 
Alternatively, let $c=q^{-1}$ in \eqref{GenAID}. Then we obtain 
\begin{tl}
    \begin{align}
        \sum_{n=0}^\infty\frac{\poq{a,b}{n}}{\poq{q}{n}\poqq{abq}{n}}q^{\frac{n(n-1)}{2}}=\frac{\poqq{a,b}{\infty}+ \poqq{aq,bq}{\infty}}{\poqq{q,abq}{\infty}}.\label{YYYYY}
    \end{align}
\end{tl}
\pf It is clear that when $c=q^{-1}$, Theorem \ref{GenAI} reduces to
\begin{align}
     \sum_{n=0}^\infty\frac{\poq{a,b}{n}q^{\frac{n(n-1)}{2}}}{\poq{q}{n}\poqq{abq}{n}}&=\frac{\poq{-1,a,b}{\infty}}{\poqq{abq}{\infty}}\sum_{n=0}^\infty\frac{b^n}{\poqq{q,-1}{n}}\sum_{k=0}^n\genfrac{[}{]}{0pt}{}{n}{k}_{q^2}(a/b)^k\nonumber\\
     &=\frac{\poq{-q,a,b}{\infty}}{\poqq{abq}{\infty}}\sum_{n=0}^\infty\frac{b^n}{\poqq{q^2}{n}}(1+q^n)\sum_{k=0}^n\genfrac{[}{]}{0pt}{}{n}{k}_{q^2}(a/b)^k\nonumber\\
    &=\frac{\poq{-q,a,b}{\infty}}{\poqq{abq}{\infty}} \big(S(a,b)+ S(aq,bq)\big),\label{LLLLL}
\end{align}
where 
\begin{align*}
    S(a,b):=\sum_{n=0}^\infty\frac{b^n}{\poqq{q^2}{n}}\sum_{k=0}^n\genfrac{[}{]}{0pt}{}{n}{k}_{q^2}(a/b)^k.
\end{align*}
It is easy to see that
\begin{align*}
    S(a,b)=\frac{1}{\poqq{a,b}{\infty}},\quad \poq{-q}{\infty}=\frac{1}{\poqq{q}{\infty}}.
\end{align*}
A substitution of these two results into \eqref{LLLLL} gives rise to \eqref{YYYYY}.\qed

With the help of \eqref{s2}, we may show Theorem \ref{GenAII} directly. 
\begin{pot2}
Let $s=0$ in Theorem \ref{andrews-s} and we appeal to Lemma \ref{lemmaadded} to evaluate $T_{n-k,k}(0)$. The result is as follows.
\begin{align*}
   \sum_{n=0}^\infty &\frac{\poq{a,b}{n}q^{\frac{n(n+1)}{2}}{c^n}}{\poq{q}{n}\poqq{ab}{n}} =\frac{\poq{a,b,-cq}{\infty}}{\poqq{ab}{\infty}}\sum_{n=0}^\infty\frac{1}{\poq{-cq}{n}}\nonumber\\
     &\qquad \times \sum_{k=0}^n\frac{\tau_2(k)b^k}{\poqq{q^2}{k}}\frac{a^{n-k}}{\poq{q}{n-2k}}    \frac{q^k+q^{n-k}}{1+q^{n}}\frac{\poq{-q^{1+n-k}}{k}}{\poq{q^{1+n-2k}}{k}}(-1)^kq^{-k^2+k}\\
     &=\frac{\poq{a,b,-cq}{\infty}}{\poqq{ab}{\infty}}\sum_{n=0}^\infty\frac{\poq{-q}{n}}{\poq{-cq}{n}}\sum_{k=0}^n\frac{b^k}{\poqq{q^2}{k}}\frac{a^{n-k}}{\poqq{q^2}{n-k}}\frac{q^k+q^{n-k}}{1+q^{n}}.
\end{align*}
Further simplification leads to
\begin{align*}
     \sum_{n=0}^\infty\frac{\poq{a,b}{n}q^{\frac{n(n+1)}{2}}{c^n}}{\poq{q}{n}\poqq{ab}{n}}&=\frac{\poq{a,b,-cq}{\infty}}{\poqq{ab}{\infty}}\sum_{n=0}^\infty\frac{1}{\poq{q,-cq}{n}}\sum_{k=0}^n\genfrac{[}{]}{0pt}{}{n}{k}_{q^2}b^ka^{n-k}\frac{q^k+q^{n-k}}{1+q^{n}}\\
     &=\frac{\poq{a,b,-cq}{\infty}}{\poqq{ab}{\infty}}\sum_{n=0}^\infty\frac{1}{\poq{q,-cq}{n}(1+q^{n})}\\
     &\qquad\times \left(\sum_{k=0}^n\genfrac{[}{]}{0pt}{}{n}{k}_{q^2}(bq)^ka^{n-k}+\sum_{k=0}^n\genfrac{[}{]}{0pt}{}{n}{k}_{q^2}b^k(aq)^{n-k}\right).
\end{align*}
Reformulate in terms of the Rogers-Szeg\"{o} polynomials \eqref{chenqi}, this yields the desired result. The theorem is proved. \qed
\end{pot2}

As a last application of Theorem  \ref{andrews-s}, we can show Theorem \ref{andrewsthree}.
\begin{pot3}
To show  \eqref{andrews-twotwo}, we first split it into two parts, as below:
\begin{align*}
    \sum_{n=0}^\infty \frac{\poq{a,b}{n}q^{\frac{n(n+1)}{2}}c^n}{\poq{q,(ab)^{1/2},-q (ab)^{1/2}}{n}}
    &=\sum_{n=0}^\infty \frac{\poq{a,b}{n}q^{\frac{n(n+1)}{2}}c^n}{\poq{q}{n}\poqq{abq^2}{n}}\frac{(1-bq^n)+(b/a)^{1/2}(1-aq^n)}{(1+(b/a)^{1/2})(1-(ab)^{1/2})}\\
    &= \frac{(1-b) L_{1}+(1-a)(b/a)^{1/2} L_{2}}{(1+(b/a)^{1/2})(1-(ab)^{1/2})}. \end{align*}
Here, we define
 \begin{align*}
   L_{1}:= \sum_{n=0}^\infty \frac{\poq{a,bq}{n}q^{\frac{n(n+1)}{2}}c^n}{\poq{q}{n}\poqq{abq^2}{n}},\quad 
  L_{2}:= \sum_{n=0}^\infty \frac{\poq{aq,b}{n}q^{\frac{n(n+1)}{2}}c^n}{\poq{q}{n}\poqq{abq^2}{n}}.
\end{align*}
From Theorem \ref{GenAI} it follows, respectively, 
   \begin{align*}
   L_{1}=
   \frac{\poq{-cq,a,bq}{\infty}}{\poqq{abq^2}{\infty}}\sum_{n=0}^\infty \frac{h_n(a,bq|q^2)}{\poq{q,-cq}{n}}, \quad L_{2}=\frac{\poq{-cq,aq,b}{\infty}}{\poqq{abq^2}{\infty}}\sum_{n=0}^\infty \frac{h_n(aq,b|q^2)}{\poq{q,-cq}{n}}.
\end{align*}
After a bit simplification, we get
 \begin{align*}
     \sum_{n=0}^\infty \frac{\poq{a,b}{n}q^{\frac{n(n+1)}{2}}c^n}{\poq{q,(ab)^{1/2},-q (ab)^{1/2}}{n}}&=\frac{1}{(1+(b/a)^{1/2})(1-(ab)^{1/2})}\frac{\poq{-cq,a,b}{\infty}}{\poqq{abq^2}{\infty}}\\
    & \times\left(
    \sum_{n=0}^\infty \frac{h_n(a,bq|q^2)}{\poq{q,-cq}{n}}+ (b/a)^{1/2}\sum_{n=0}^\infty \frac{h_n(aq,b|q^2)}{\poq{q,-cq}{n}}\right).
\end{align*}
Dividing both sides by  $1+(ab)^{1/2}$ leads us to the desired identity.
\qed
\end{pot3}

The case $c=1$ of Theorem \ref{andrewsthree} yields the following result which is very similar to \cite[Thm. 1]{222}. 
\begin{tl}\label{andrews++}
\begin{align}
    \sum_{n=0}^\infty \frac{\poq{a,b}{n}q^{\frac{n(n+1)}{2}}}{\poq{q,(ab)^{1/2}}{n}\poq{-(ab)^{1/2}}{n+1}}&=\frac{a^{1/2}\poqq{aq,b}{\infty}+b^{1/2}\poqq{a,bq}{\infty}}{(a^{1/2}+b^{1/2})\poqq{q,ab}{\infty}}.\label{andrews-twotwotwo}
\end{align}
\end{tl}
\pf When $c=1$, the conclusion follows from \eqref{andrews-twotwo} directly by using \eqref{newGF}. 
\qed

\section{New proof of Warnaar's partial theta identity}\label{sec3}
Recall that in their paper \cite{Warnaarand}, Andrews and Warnaar established the following beautiful partial theta identity.
\begin{yl}[Cf. {\rm \cite[Thm. 1.5]{Warnaarhim}}]
    \begin{align}
1+\sum_{n=1}^{\infty}(-1)^n q^{\binom{n}{2}}\left(a^n+b^n\right)=(q, a, b ; q)_{\infty} \sum_{n=0}^{\infty} \frac{(a b / q ; q)_{2 n}}{(q, a, b, a b ; q)_n} q^n.\label{MMMMMM-1}
\end{align}
\end{yl}
It is worth mentioning that in their paper \cite{wangma}, Wang and Ma put forward a unified method to such kind of theta identities.  In this part, we will focus on the special case of \eqref{MMMMMM-1} and  treat it via the foregoing argument. 
\begin{dl}[Warnaar's partial theta identity: {\rm \cite[p. 4]{Warnaarhim}}]\label{war-one} 
    \begin{align}
 1+2 \sum_{n=1}^{\infty} a^n q^{2 n^2}=(q ; q)_{\infty}\poqq{aq}{\infty}\sum_{n=0}^{\infty} \frac{(-a ; q)_{2 n} q^n}{(q,-a q ; q)_n\poqq{aq}{n}}.    \label{wwwwww}
\end{align}
\end{dl}
During this procedure, we come to a new Bailey pair and some allied $q$-identities lying behind this theta identity. To this end, we recall first  the well known Bailey lemma  associated with the Bailey pairs.
\begin{yl}[Bailey lemma: \mbox{\cite[Chap. 3]{111-111}}]\label{baileyyl}
	For any integer $n\geq 0$, it holds
	\begin{align}
		\frac{1}{\left(a q / \rho_1, a q / \rho_2 ; q\right)_n} \sum_{k=0}^n \frac{\left(\rho_1, \rho_2 ; q\right)_k\left(a q / \rho_1 \rho_2 ; q\right)_{n-k}}{(q ; q)_{n-k}}\left(\frac{a q}{\rho_1 \rho_2}\right)^k \beta_k(a, q) \nonumber\\
		=\sum_{k=0}^n \frac{\left(\rho_1, \rho_2 ; q\right)_k}{(q ; q)_{n-k}(a q ; q)_{n+k}\left(a q / \rho_1, a q / \rho_2 ; q\right)_k}\left(\frac{a q}{\rho_1 \rho_2}\right)^k \alpha_k(a, q),\label{bailey}
	\end{align}
where $(\alpha_n(a,q),\beta_n(a,q))$ is a Bailey pair with respect  to $a$ and $q$, defined by
	\begin{align}
	\beta_n(a,q)=\sum_{k=0}^{n}\frac{\alpha_k(a,q)}{(q;q)_{n-k}(aq;q)_{n+k}}.\label{baileybeta}
	\end{align}
\end{yl}
Using  this lemma, we can show
\begin{xz} Define, for any integer $n\geq 0$, the finite sum
    \begin{align}
   \gamma(n):= \sum_{k=0}^n\genfrac{[}{]}{0pt}{}{n}{k}_q\frac{ q^{k^2}}{\poq{-q}{k}}.\label{finalfinal}
     \end{align}
   Then 
\begin{align}(\alpha_n(1,q),\beta_n(1,q))=\bigg(2(-1)^nq^{2n^2},\frac{1}{\poq{q}{n}^2}+ \frac{\gamma(n)}{\poq{q}{n}}\bigg)\label{great111-321}
\end{align}
is a Bailey pair with respect to $a=1$.
\end{xz}
\pf First of all, letting
$\rho_1,\rho_2\to \infty$ and $a=1$, then \eqref{bailey} reduces to
\begin{align}
    \sum_{k=0}^n \frac{q^{k^2}}{\poq{q}{n-k}} \beta_{k}(1,q)=\sum_{k=0}^n\frac{q^{k^2}}{\poq{q}{n+k}\poq{q}{n-k} }\alpha_{k}(1,q).\label{baileysp}
\end{align}
In such case, we see that \eqref{baileybeta} reduces to
\begin{align}
    \beta_n(1,q)=\sum_{k=0}^{n}\frac{\alpha_k(1,q)}{(q;q)_{n-k}(q;q)_{n+k}}.\label{pairsp}
\end{align}
In order to show that \eqref{great111-321} is a Bailey pair, we first set
\begin{align*}
    \beta_n(1,q):=\frac{1}{\poqq{q^2}{n}}+\frac{1}{\poq{q}{n}^2}.
\end{align*}
By the inverse relation,  we may solve   $\alpha_n(1,q)$ from \eqref{pairsp}, as follows
\begin{align*}
    \frac{\alpha_n(1,q)}{(1+q^n)\tau(n)}=\sum_{k=0}^{n}\frac{\poq{q^{-n},q^n}{k}}{\poq{q}{k}^2}q^k+\sum_{k=0}^{n}\frac{\poq{q^{-n},q^n}{k}}{\poqq{q^2}{k}}q^k,
\end{align*}
which, after an application of the $q$-Chu-Vandermonde formula \eqref{s1}, turns  out to be
\begin{align*}
    \alpha_n(1;q)&=(1+q^n)\tau(n)\frac{\poq{-q^{1-n}}{n}}{\poq{-q}{n}}q^{n^2}=2(-1)^nq^{n^2}.
\end{align*}
Therefore, we come to a Bailey pair 
\begin{align}
\left(2(-1)^nq^{n^2},\frac{1}{\poqq{q^2}{n}} + \frac{1}{\poq{q}{n}^2}\right).\label{3666Bailey}
\end{align} 
A direct substitution of \eqref{3666Bailey} into  \eqref{baileysp} gives rise to \eqref{great111-321}.
\qed

With the above setup,  we are now in good position to show Warnaar's theta identity \eqref{MMMMMM-1}.

\begin{pot5} It suffices to evaluate the sum on the right-hand side of \eqref{wwwwww}. At first,  we see
\begin{align*}
   \mbox{RHS of \eqref{wwwwww}}&= (q ; q)_{\infty}\left(a q ; q^2\right)_{\infty} \sum_{n=0}^{\infty} \frac{(-a ; q)_{2 n} q^n}{(q,-a q ; q)_n\left(a q ; q^2\right)_n}\\
    &=\frac{\poq{q,-a}{\infty}}{\poq{-aq}{\infty}}  \sum_{n=0}^{\infty}\frac{q^n}{\poq{q}{n}}\frac{\poqq{aq^{1+2n}}{\infty}}{\poq{-aq^{2n}}{\infty}}\poq{-aq^{1+n}}{\infty}.
    \end{align*}
    By  Lemma \ref{main}, we have
    \begin{align*}
    \text{RHS of \eqref{wwwwww}}
    &=\poq{q}{\infty}(1+a) \sum_{n,M=0}^\infty\frac{q^n}{\poq{q}{n}} \sum_{i+j=M}\frac{q^{(1+n)i}q^{\binom{i}{2}}}{\poq{q}{i}}\lambda_j(-q)(-q^{2n})^ja^M\\
  & =\poq{q}{\infty}(1+a) \sum_{M\geq j \geq 0} \frac{q^{\binom{M-j}{2}}q^{M-j}}{\poq{q}{M-j}}(-1)^j\lambda_j(-q)a^M\sum_{n=0}^{\infty}\frac{q^{(1+M+j)n}}{\poq{q}{n}}\\
    &=(1+a) \sum_{M\geq j \geq 0} \frac{\poq{q}{M+j}}{\poq{q}{M-j}}(-1)^{j}q^{\binom{M-j}{2}}q^{M-j}\lambda_j(-q)a^M,
\end{align*}
where $\lambda_j(-q)$ is given by \eqref{coeff}. Relabeling $M$ by $n$ and $j$ by $k$, we get
\begin{align*}
    \mbox{RHS of \eqref{wwwwww}}
    &=(1+a) \sum_{n=0}^\infty a^nq^{n+\binom{n}{2}}\sum_{k=0}^n\frac{\poq{q}{n+k}}{\poq{q}{n-k}} (-1)^{k}q^{\binom{k}{2}-nk}\lambda_k(-q)\\
    &=\sum_{n=0}^\infty a^nF(n)+\sum_{n=0}^\infty a^{n+1}F(n)=1+\sum_{n=1}^\infty a^n(F(n)+F(n-1)),
\end{align*}
where
\begin{align*}
    F(n):=q^{n+\binom{n}{2}}\sum_{k=0}^n\frac{\poq{q}{n+k}}{\poq{q}{n-k}}(-1)^kq^{\binom{k}{2}-nk}\lambda_k(-q).
\end{align*}
    It is not hard to check that
\begin{align*}
    F(n)+F(n-1)
    &=q^{\binom{n}{2}}\sum_{k=0}^{n}\frac{\poq{q}{n+k-1}}{\poq{q}{n-k}}\tau(k)q^{-nk}\lambda_{k}(-q)\big\{q^n(1-q^{n+k})+q^k(1-q^{n-k})\big\}\\
    &=q^{\binom{n}{2}}(1-q^{2n})\sum_{k=0}^{n}\frac{\poq{q}{n+k-1}}{\poq{q}{n-k}}\tau(k)q^{k-nk}\lambda_{k}(-q)\\
       &=q^{\binom{n}{2}}(1+q^n)\sum_{k=0}^{n}\poq{q^{-n},q^n}{k}q^k\lambda_{k}(-q).
\end{align*}
On comparing with the left-hand side of \eqref{wwwwww}, clearly we only need to show 
\begin{align}
   q^{\binom{n}{2}}(1+q^n)\sum_{k=0}^{n}\poq{q^{-n},q^n}{k}q^k\lambda_{k}(-q)=2q^{2n^2}.\label{success}
\end{align}
The argument goes as follows.  First, we recast \eqref{success} as the form 
\begin{align}
    \beta_n=\sum_{k=0}^{n}\poq{q^{-n},q^n}{k}q^k\lambda_{k}(-q),\label{oldform}
\end{align}
where
\begin{align*}
    \left\{\begin{aligned}
          &\lambda_{n}(-q)=\displaystyle\frac{q^{n^2}}{\poqq{q^2}{n}}{}_2\phi_1\left(\begin{array}{cc}
        q^{-n}, & -q^{-n} \\
         & 0
    \end{array};q,-q\right) \\
          &\beta_n= \frac{2q^{2n^2-\binom{n}{2}}}{1+q^n}\quad(n\geq 1); \quad \beta_0=1.
    \end{aligned}\right.
\end{align*}
By virtue of  Carlitz's matrix inversion
\begin{align*}
\bigg(\frac{\left(q^{-n}, a q^n ; q\right)_k}{(q, a q ; q)_k} q^k\bigg)_{n\geq k\geq 0}^{-1}=\bigg(\frac{\left(a, q^{-n} ; q\right)_k}{\left(q, a q^{1+n} ; q\right)_k} \frac{1-a q^{2 k}}{1-a} q^{k n}\bigg)_{n\geq k\geq 0},
\end{align*}
it is easily seen that \eqref{oldform} is equivalent to
\begin{align*}
    \poq{q}{n}^2\lambda_{n}(-q)&=1+\sum_{k=1}^n \frac{\poq{q^{-n}}{k}\poq{q}{k-1}}{\left(q, q^{1+n} ; q\right)_k}(1-q^{2 k}) q^{k n}\beta_k.
\end{align*}
That is
\begin{align*}
 \lambda_{n}(-q)&=\frac{1}{ \poq{q}{n}^2}+\frac{1}{ \poq{q}{n}}\sum_{k=1}^n \frac{\poq{q^{-n}}{k}\poq{q}{k-1}}{\left(q;q)_k(q;q\right)_{n+k}}(1-q^{2 k}) q^{k n}\beta_k
 \\
     &=\frac{1}{ \poq{q}{n}^2}+\sum_{k=1}^n \frac{\tau(k)(1+q^k)}{(q;q)_{n-k}(q;q)_{n+k}}\beta_k=\sum_{k=-n}^n \frac{(-1)^kq^{2k^2}}{(q;q)_{n-k}(q;q)_{n+k}}.
\end{align*}
Lastly, we get
\begin{align*}
    \frac{q^{n^2}}{\poqq{q^2}{n}}\sum_{k=0}^n\frac{\poq{q^{-n},-q^{-n}}{k}}{\poq{q}{k}}(-q)^k=\sum_{k=-n}^n \frac{(-1)^kq^{2k^2}}{(q;q)_{n-k}(q;q)_{n+k}}.
\end{align*}
Using the basic relation \cite[(I. 10)]{10}
\[(-q^{-n};q)_k=\frac{\poq{-q}{n}}{\poq{-q}{n-k}}(-1)^kq^{-nk}\tau(k)\]
to simplify the last identity, only finding that it is the Bailey lemma given by the Bailey pair \eqref{3666Bailey}.
Summing up, \eqref{success} is proved.\qed
\end{pot5}

We end our discussion by the following $q$-identities via combination of two Bailey pairs  \eqref{great111-321} and \eqref{3666Bailey} with the Bailey lemma. 

\begin{xz}With the same notation as Lemma \ref{baileyyl} and $\gamma(n)$ be given by \eqref{finalfinal}. Then we have
\begin{align}
		\frac{1}{\left(q / \rho_1,  q / \rho_2 ; q\right)_n} \sum_{k=0}^n \frac{\left(\rho_1, \rho_2 ; q\right)_k\left(  q / \rho_1 \rho_2 ; q\right)_{n-k}}{\poq{q}{k}(q ; q)_{n-k}}\left(\frac{  q}{\rho_1 \rho_2}\right)^k \bigg(\frac{1}{\poq{q}{k}}+ \gamma(k)\bigg) \nonumber\\
		=2\sum_{k=0}^n \frac{\left(\rho_1, \rho_2 ; q\right)_kq^{2k^2}}{(q ; q)_{n-k}(  q ; q)_{n+k}\left(  q / \rho_1,   q / \rho_2 ; q\right)_k}\left(\frac{ - q}{\rho_1 \rho_2}\right)^k.\label{bailey-concrete1}
	\end{align}
    In particular
\begin{description}
   \item [(i) {\rm ($ \rho_1=1/a,  \rho_2=q$ in \eqref{bailey-concrete1})}] For $|a|<1$, it holds 
   \begin{align*}
		\sum_{k=0}^n \frac{\left(1/a; q\right)_k\left( a; q\right)_{n-k}}{(q ; q)_{n-k}}a^{k}\gamma(k)\\
=2\poq{q}{n-1}\left(aq; q\right)_{n}\sum_{k=0}^n \frac{\left(1/a; q\right)_k(1-q^k)q^{2k^2}}{\left(aq; q\right)_{k}(q ; q)_{n-k}(  q ; q)_{n+k}}(-a)^{k}.\nonumber
	\end{align*}
    The limitation $n\to \infty$ yields
    \begin{align*}
		(1-a)\sum_{k=0}^\infty\left(1/a; q\right)_ka^{k}\gamma(k)=2 \sum_{k=0}^\infty \frac{\left(1/a; q\right)_k}{\left(aq; q\right)_{k}}(-a)^{k}(1-q^k)q^{2k^2}.
	\end{align*}
    
 \item[(ii) {\rm ($ \rho_1,\rho_2\to \infty$ in \eqref{bailey-concrete1})}] \begin{align*}
		 \sum_{k=0}^n \frac{q^{k^2}\gamma(k)}{\poq{q}{k}(q ; q)_{n-k}}
		=\sum_{k=-n}^n \frac{(-1)^kq^{3k^2}}{(q ; q)_{n-k}(q ; q)_{n+k}},
	\end{align*}
    whose limitation $n\to \infty$ yields
    \begin{align}
		 \sum_{k=0}^\infty \frac{q^{k^2}\gamma(k)}{\poq{q}{k}}
		=\frac{(q^3;q^6)_\infty}{(q,q^2;q^3)_\infty}.
        \label{bailey-concrete222}
	\end{align}
\end{description}
 \end{xz}
 \pf Actually, \eqref{bailey-concrete1} is the direct consequence of the Bailey lemma specialized by the Bailey pair \eqref{great111-321}, while \eqref{bailey-concrete222} involves the Jacobi triple product identity
 \[\sum_{k=-\infty}^\infty (-1)^kq^{3k^2}=(q^3,q^3,q^6;q^6)_\infty.\]
  \qed

When the Bailey pair \eqref{3666Bailey} is taken into account, we also have  
\begin{xz}With the same notation as Lemma \ref{baileyyl}. We have
\begin{align}
		\frac{1}{\left(q / \rho_1,  q / \rho_2 ; q\right)_n}\sum_{k=0}^n \frac{\left(\rho_1, \rho_2 ; q\right)_k\left(  q / \rho_1 \rho_2 ; q\right)_{n-k}}{(q ; q)_{n-k}}\left(\frac{  q}{\rho_1 \rho_2}\right)^k\frac{(q ; q)_{k}+(-q ; q)_{k}}{\poqq{q^2}{k}(q ; q)_{k}}\nonumber\\
		=2\sum_{k=0}^n \frac{\left(\rho_1, \rho_2 ; q\right)_kq^{k^2}}{(q ; q)_{n-k}(  q ; q)_{n+k}\left(  q / \rho_1,   q / \rho_2 ; q\right)_k}\left(\frac{ - q}{\rho_1 \rho_2}\right)^k.\label{bailey-concrete0-123}
	\end{align}
   In particular
\begin{description}
   \item [(i) {\rm ($ \rho_1=1/a,  \rho_2=q$  in \eqref{bailey-concrete0-123})}] For $|a|<1$, it holds \begin{align*}
		\sum_{k=0}^n \frac{\left(1/a; q\right)_k\left(a; q\right)_{n-k}}{(q ; q)_{n-k}\poq{-q}{k}}a^k\\
=2\poq{q}{n-1}\left( aq; q\right)_{n}\sum_{k=0}^n \frac{\left(1/a; q\right)_k(1-q^k)q^{k^2}}{\left( aq; q\right)_{k}(q ; q)_{n-k}( q ; q)_{n+k}}\left(-a\right)^{k}.\nonumber
	\end{align*}
 \item[(ii) {\rm ($ \rho_1,\rho_2\to \infty$ in \eqref{bailey-concrete0-123})}] \begin{align}
		\sum_{k=0}^n \frac{q^{k^2}}{\poqq{q^2}{k}(q ; q)_{n-k}}=\sum_{k=-n}^n \frac{(-1)^kq^{2k^2}}{(q ; q)_{n-k}(  q ; q)_{n+k}}.\label{bailey-concrete000}
	\end{align}
\end{description}
\end{xz}
\pf Clearly, \eqref{bailey-concrete0-123} is the direct consequence of the Bailey lemma specialized by the Bailey pair \eqref{3666Bailey}.  Note that \eqref{bailey-concrete000} comes from
 \[ \sum_{k=0}^n \frac{q^{k^2}}{\poq{q}{k}^2(q ; q)_{n-k}}
		=\frac{1}{\poq{q}{n}^2}.\]
 \qed

 \noindent{\bf Funding} \quad This research was supported by the Natural Science Foundation of Zhejiang Province, Grant No. LY24A010012 and the National Natural Science Foundation of China, Grant No.  12471315.
 
\noindent{\bf Declarations Conflict of interest}\quad The author declares that they have no competing interests related to this work.
\bibliographystyle{amsplain}

\end{document}